\documentclass[12pt]{amsart}
\usepackage{latexsym, amsmath,amssymb}
\usepackage{hyperref}
\usepackage{xcolor}
\hypersetup{
    colorlinks,
    linkcolor={red!50!black},
    citecolor={blue!50!black},
    urlcolor={blue!80!black}
}

\usepackage{graphicx}

 \numberwithin{equation}{section}

\def\XXint#1#2#3{{\setbox0=\hbox{$#1{#2#3}{%
\int}$ }
\vcenter{\hbox{$#2#3$ }}\kern-.6\wd0}}

\renewcommand{\epsilon}{\varepsilon}
\newtheorem{theorem}{Theorem}

\newtheorem{lemma}[theorem]{Lemma}
\newtheorem{corr}[theorem]{Corollary}

\newtheorem{proposition}[theorem]{Proposition}
\newtheorem{deff}[theorem]{Definition}
\newtheorem{remark}[theorem]{Remark}
\newcommand{\bth}{\begin{theorem}}
\newcommand{\ble}{\begin{lemma}}
\newcommand{\bcor}{\begin{corr}}

\newcommand{\bdeff}{\begin{deff}}

\newcommand{\bprop}{\begin{proposition}}
\newcommand{\ele}{\end{lemma}}
\newcommand{\ecor}{\end{corr}}
\newcommand{\edeff}{\end{deff}}

\numberwithin{theorem}{section}

\newcommand{\eprop}{\end{proposition}}

\renewcommand{\Pi}{\varPi}

\renewcommand{\epsilon}{\varepsilon}

\begin{document}

\title[Integral representation for fractional Laplace operators]
{Integral representation for\\ fractional Laplace-Beltrami operators}
\author[D. Alonso-Or\'an]{Diego Alonso-Or\'an}
\address{Instituto de Ciencias Matem\'aticas CSIC-UAM-UC3M-UCM -- Departamento de Matem\'aticas (Universidad Aut\'onoma de Madrid), 28049 Madrid, Spain} 
\email{diego.alonso@icmat.es}
\author[A. C\'ordoba]{Antonio C\'ordoba}
\address{Instituto de Ciencias Matem\'aticas CSIC-UAM-UC3M-UCM -- Departamento de Matem\'aticas (Universidad Aut\'onoma de Madrid), 28049 Madrid, Spain} 
\email{antonio.cordoba@uam.es}
\author[A. D. Mart\'inez]{\'Angel D. Mart\'inez}
\address{Instituto de Ciencias Matem\'aticas (CSIC-UAM-UC3M-UCM) -- Departamento de Matem\'aticas (Universidad Aut\'onoma de Madrid), 28049 Madrid, Spain} 
\email{angel.martinez@icmat.es}

\begin{abstract}
In this paper we provide an integral representation of the fractional Laplace-Beltrami operator for general riemannian manifolds which has several interesting applications. We give two different proofs, in two different scenarios, of essentially the same result. One of them deals with compact manifolds with or without boundary, while the other approach treats the case of riemannian manifolds without boundary whose Ricci curvature is uniformly bounded below.
\end{abstract}

\maketitle

\section{\textbf{Introduction}}

Fractional laplacians appear in a number of equations of mathematical interest. In the euclidean space and tori explicit expresions are well-known. For example Zygmund's operator, $\sqrt{-\Delta}$, can be expressed in the tori as a principal value integral
\[\Lambda f(x)=c_n P.V.\sum_{\nu\in\mathbb{Z}^n}\int_{\mathbb{T}^n}\frac{f(x)-f(y)}{|x-y-\nu|^{n+1}}dy\]
while in the euclidean space we have
\[\Lambda f(x)=c_n P.V.\int_{\mathbb{R}^n}\frac{f(x)-f(y)}{|x-y|^{n+1}}dy.\]
This expresions combined with elementary algebra provide useful pointwise estimates (see, for instance, \cite{CC, CC2, CV}). In other situations lack of explicit expresions like the one stated above make the analysis much harder to achieve, this typically corresponds to situations involving boundary effects or anisotropy.

In this paper we will provide an integral representation of the fractional Laplace-Beltrami operator on a general compact manifold with a nice error term. As a first example of the power of these explicit formulae we present direct proofs of fractional Sobolev's embeddings on compact manifolds. Furthermore, the explicit kernels make transparent the r\^ole played by different ad hoc definitions of fractional integration. The authors have applied these representations in their recent work \cite{AOCM2} to prove global existence of the critical surface quasigeostrophic equation on the two dimensional standard sphere. The present work stems from our original approach \cite{CM} which remained unpublished since we found another alternative proof which we judged to be specially elegant. But one of the features of the present approach is that it allows to improve the C\'ordoba-C\'ordoba inequality (cf. \cite{CC, CC2}) in the spirit of the nonlinear lower bounds due to Constantin and Vicol (cf. \cite{CVi}), namely, one may achieve pointwise estimates of the form
\[\nabla f(x)\cdot\Lambda^{\alpha}\nabla f(x)\geq\frac{1}{2}|\nabla f(x)|^2+\frac{|\nabla f(x)|^{2+\alpha}}{c\|f\|^{\alpha}_{L^{\infty}(\mathbb{R}^n)}}.\]
Those turn out to be quite useful concerning the important question of global existence of solutions to the critical surface quasigeostrophic equation in $\mathbb{R}^n$. Their work followed landmark results obtained independently by Kiselev, Nazarov and Volberg \cite{KNV} and Caffarelli and Vasseur \cite{CV}. We believe that its usefulness justifies now its publication due to its instrumental r\^ole in the subsequent contributions of the authors to the regularity of weak solutions to the (critical) drift diffusion equations on manifolds, namely
\[\partial_t\theta+u\cdot\nabla_g\theta=-\Lambda\theta\]
where $u$ is a divergence free vector field that depends on $\theta$ in some specific manner (e.g. Riesz's transform). Let us highlight the following result for the critical surface quasigeostrophic equation on the sphere

\begin{theorem}\label{sphere}
Let the vector field $u=\nabla_g^{\perp}\Lambda^{-1}\theta$ and the initial data $\theta_0\in C^{\infty}(\mathbb{S}^2)$. Then the unique global weak solution is smooth for all times, i.e. $\theta\in C^{\infty}([0,\infty)\times\mathbb{S}^2)$.
\end{theorem}
Its proof will be published elsewhere. It combines De Giorgi's ideas as in the work of Caffarelli and Vasseur with some nonlinear maximum principles coming from the work of Constantin and Vicol (cf. \cite{CV, CVi}). Let us add that in some related work of Constantin and Ignatova concerning the equation on euclidean bounded domains \cite{CI1, CI2} they regret the lack of explicit expresions for their kernels. Our method can be used in their context, but allowing as well to include the difficulty produced by anisotropy.

\subsection{The integral representation}

As usual let $(M,g)$ be a closed compact manifold of dimension $n\geq 2$ whose Laplace-Beltrami operator is denoted by $-\Delta_g$. The following is the main result of this paper

\begin{theorem}\label{prop}
Let $f$ be smooth and $s\in(0,1)$, then for any fixed $N$ one has the representation
\[(-\Delta_g)^{s}f(x)=\textrm{P.V.}\int_M\frac{f(x)-f(y)}{d(x,y)^{n+2s}}(c_{n,s}\chi u_0+k_N)(x,y)d\textrm{vol}_g(y)+O(\|f\|_{H^{-N}(M)}),\]%+\int_M\frac{g(y)-g(x)}{d(x,y)^{n-1+2s}}k(x,y)d\textrm{vol}_g(y)\]
where $k_N(x,y)=O(d(x,y))$ is a smooth function, $\chi$ is a smooth cut off function equal to one on the diagonal, the implicit constant depends on $N$, $c_{s,n}>0$ is a constant independent of $N$ and $u_0$ is a smooth function such that $u_0(x,x)=1$.%where $g=f+Sf$ for some operator $S$ that gains derivatives and $k$ is a certain continuous function. Both of them which might be expressed rather explicitly.
\end{theorem}

Notice that the norm in the error might be taken to be $L^{\infty}$. The proof uses spectral calculus interwined with Hadamard parametrix to provide an explicit integral representation, similar to a classical singular integral, with a harmless error. This is enough for the applications we have in mind and in order to establish it we will start with the following

\begin{lemma}
For $s\in (-1,0)$ and under the hypotheses of the previous theorem
\[(-\Delta_g)^{s}f(x)=\int_M\frac{f(y)}{d(x,y)^{n+2s}}(c_{n,s}\chi u_0+k_N)(x,y)d\textrm{vol}_g(y)+O(\|f\|_{H^{-N}(M)})\]
holds.
\end{lemma}
An integration by parts argument will then imply the theorem for any $s\in (0,1)$. Notice that in the two dimensional setting for $s=-1$ a logarithm singularity must appear in the lemma, but here we will not consider this case which need such well known modification. The explicit expression of the error is crucial for the applications, and one may compare it with similar expresions for the cases of the flat tori or euclidean space, as was pointed out in the introduction (cf. \cite{CC,CC2,S1}). In the case of manifolds with boundary there is controlled degeneration of constants as one approaches the boundary. 

In the non compact case, instead of Hadamard parametrix we will use the fundamental solution of the heat operator which, nevertheless, is related to the former by a Laplace transform (cf. \cite{MP}). One needs then to employ a different identity to begin with, together with some sharp heat kernel bounds whose validity demands geometrical restrictions on the Ricci curvature, which has to be bounded below. Also a strictly positive lower bound for the injectivity radius is required. 

\begin{theorem}\label{prop2}
Let $(M,g)$ be a closed riemannian manifold with Ricci curvature and injectivity radius bounded below, $f$ be smooth and $s\in(0,1)$. Then one has the following representation
\[(-\Delta_g)^{s}f(x)=\textrm{P.V.}\int_M\frac{f(x)-f(y)}{d(x,y)^{n+2s}}(c_{n,s}\chi U_0+k)(x,y)d\textrm{vol}_g(y)+O(\|f\|_{\infty}),\]%+\int_M\frac{g(y)-g(x)}{d(x,y)^{n-1+2s}}k(x,y)d\textrm{vol}_g(y)\]
where $k(x,y)=O(d(x,y))$ is a smooth function, $\chi$ is a smooth diagonal cut off function, $c_{s,n}>0$ and $U_0$ is a smooth function such that $U_0(x,x)=1$.
\end{theorem}

%In the proof the error term can be taken to be $O(\|f\|_2)$ but it does not seem to be possible to take a smoother error term. 

Observe that the error term has a remarkable difference with the main result stated before (Theorem \ref{prop}). We will present this result first, which is more general although not quite as strong as the former. Afterwards we will introduce the needed notational convention together with the essentials of the Hadamard parametrix construction. Then we will end the paper with a proof of the main result and, for the sake of completeness, the fractional Sobolev embedding theorem that was needed in the proof of Theorem \ref{sphere}.

\begin{corr}
Let $(M,g)$ be a compact riemannian manifold of dimension $n$, $s\in (0,\frac{1}{2})$ and $p=\frac{2n}{n-2s}$. Then there exist a constant $C>0$ such that
\[\|f\|_p\leq C(\|f\|_2+\|\Lambda^sf\|_2).\]
\end{corr}

The cases when $s$ is an integer are well known and can be found in \cite{Au, Au2}.

\section{\textbf{Proof of Theorem \ref{prop2}: the non compact case}}

In this proof we will take advantage of the following  well known representation formula
\[(-\Delta_g)^{s}f(x)=\int_0^{\infty}(e^{-t\Delta_g}f(x)-f(x))\frac{dt}{t^{1+s}}.\]
The semigroup action might be expressed through the heat kernel $G(x,y,t)$ as follows
\[\int_0^{\infty}\int_MG(x,y,t)(f(y)-f(x))d\textrm{vol}_g(y)\frac{dt}{t^{1+s}}.\]
We will split this integral in three parts: corresponding to large times, small times but far from the spatial singularity and, finally, small times near the singularity. In fact the latter is the main part contributing to the kernel in our expression while the other two will add to the error term in the statement. The first integral corresponds to
\[\int_1^{\infty}\int_MG(x,y,t)(f(y)-f(x))d\textrm{vol}_g(y)\frac{dt}{t^{1+s}}\]
which is $O(\|f\|_{\infty})$. The second integral has the form
\[\int_0^{1}\int_{d_g(x,y)>1}G(x,y,t)(f(y)-f(x))d\textrm{vol}_g(y)\frac{dt}{t^{1+s}}.\]
To control it we use the bound of Li and Yau (Corollary 3.1 from \cite{LY}). Here the condition on the curvature arises, namely, their bound assures the existence of a positive constant $C$ such that
\[G(x,y,t)\leq C \frac{1}{\textrm{vol}^{1/2}_g(B_x(\sqrt{t}))\textrm{vol}^{1/2}_g(B_y(\sqrt{t}))}e^{-C\kappa t}e^{-\frac{d_g(x,y)^2}{5t}}\]
where $\textrm{Ric}_g\geq \kappa$, $C$ is some positive constant and $B_x(r)$ denotes the ball of radius $r$ centered at $x$. It is quite elementary to check that such an integral is bounded by $O(\|f\|_{\infty})$. We are hence left with
\[\lim_{\epsilon\rightarrow 0}\int_0^{1}\int_{\epsilon<d_g(x,y)<1}G(x,y,t)(f(y)-f(x))d\textrm{vol}_g(y)\frac{dt}{t^{1+s}}.\]
In such a range we can use the heat kernel parametrix (which is actually closely related to Hadamard's parametrix, cf. \cite{MP}), namely:
\[G(x,y,t)=\frac{1}{(4\pi t)^{n/2}}e^{-\frac{d_g(x,y)^2}{4t}}\left(U_0(x,y)\chi(x,y)+\chi(x,y)\sum_{j=1}^k U_j(x,y)t^j+O(t^{k+1})\right),\]
where $U_0(x,x)=1$, $\chi$ is a radial cut off function around $x$ and $k$ is big enough. We refer to \cite{MP} or \cite{Ch} for further details. One may expand the $O$-term asymtotically where further powers of $t$ will appear together with some smooth functions having an interesting geometric meaning. However, they are unnecessary for our current purposes. Plugging that information into the above integral one finds after a trivial change of variables that it is equal to:
\[\lim_{\epsilon\rightarrow 0}\int_{\epsilon<d_g(x,y)<1}\frac{f(y)-f(x)}{d(x,y)^{n+\alpha}}\int_0^{1/d(x,y)^2}\frac{e^{-1/t}(\chi(x,y)U_0(x,y)+O(t))}{t^{n/2+1+2s}}dt d\textrm{vol}_g(y).\]
One can now complete the range of the inner integral to the whole interval $(0,\infty)$ to obtain a constant, while noticing that the added part decays very rapidly to zero. %Finally, to get smooth properties for the kernel one just needs to perform the above spliting the integrals using some smooth cut off.

\section{\textbf{The Hadamard parametrix}}

In this section we present a brief description of the Hadamard parametrix construction which is suited for our purposes. A complete though rather technical treatment can be found in H\"ormander's \cite{Ho}. A more recent and accesible exposition can be found in \cite{DKS}. We intermingle both here. Let us introduce $F_0^z$ known as a Bessel potential of order zero which is a fundamental solution of $(-\Delta-z)F_0^z(s)=\delta_0(x)$. Fourier analysis can be employed to obtain the following representation formula
\[F_0^z(x)=(2\pi)^{-n}\int_{\mathbb{R}^n}e^{ix\cdot\xi}(|\xi|^2-z)^{-1}d\xi.\]
Notice that it is radial being the Fourier transform of another radial function. The Hadamard parametrix method introduces more potentials, $F_{\nu}^z(x)$, as part of the construction which can be expressed in terms of modified Bessel functions of the second kind as follows (cf. \cite{DKS})
\[F_{\nu}^z(x)=c_{\nu}|x|^{-\frac{n}{2}+\nu+1}z^{\frac{n}{4}-\frac{\nu+1}{2}}K_{n/2-\nu-1}(\sqrt{z}|x|).\]
The special functions involved satisfy the following bound for $\textrm{Re}(w)>0$:
\[|K_{\ell}(w)|\leq\left\{\begin{array}{cl}
Cw^{-\ell}&\textrm{if $|w|\leq 1$}\\
Ce^{-w}&\textrm{if $|w|>1$}
\end{array}\right.\]
for some absolute constant $C>0$ (cf. \cite{W}, $\S7\cdot 23$). They also satisfy several recursive relations among which we will use $-2\frac{\partial F_{\nu}}{\partial x}=x F_{\nu-1}$ for $\nu>0$. Let $u_0$ be some function to be specified later
\[(-\Delta_g-z)(u_0F_0^z)=u_0(0)\textrm{det}(g^{ik})^{\frac{1}{2}}\delta_0+(-\Delta_gu_0)F_0+2(hu_0-2\langle x,\frac{\partial u_0}{\partial x}\rangle)(F_0^z)'(|x|^2).\]
where $h(x)=\sum g_{jk}(x)b^j(x)x_k$ in normal coordinates. Here we are employing normal coordinates and a consequence of Gauss lemma implicitly. To get rid of the last term we will know choose $u_0$ to be a function such that $u_0(0)=1$ and
\[hu_0=2\langle x,\frac{\partial u_0}{\partial x}\rangle.\]
For further approximations one proceeds similarly, where all the functions $u_{\nu}$ that appear in the process are smooth. We refer the reader to the aforementioned references for further details on the construction. Let $\chi(x,y)$ be a cut off function supported near the diagonal, it follows that the operator
\[\mathcal{P}_N^zf(x)=\int_M\chi(x,y)\sum_{\nu=0}^Nu_{\nu}(x,y)F_{\nu}^z(d(x,y))f(y)d\textrm{vol}_g(y)\]
is a right parametrix of $(-\Delta_g-z)$ that is $(-\Delta_g-z)\mathcal{P}^z_N=\delta+R^z$ where
\[R_N^zf(x)=\int_M \chi(x,y)h_N(x,y)F^z_{N}(d_g(x,y))f(y)d\textrm{vol}_g(y)\]
where $h$ is some specific smooth function.

\section{\textbf{Proof of Theorem \ref{prop}: the compact case}}

Let $s<0$. We may use spectral calculus to define
\[(-\Delta_g)^{s}f(x)=\frac{1}{2\pi i}\int_{\gamma}z^{s}(-\Delta_g-z)^{-1}f(x)dz\]
where $\gamma$ is some appropiate contour which avoids the negative real numbers, i.e. we choose a branch of the logarithm so that $z^{\alpha}$ is holomorphic. This countour integral must be interpreted as a principal value integral on which $\gamma$ corresponds to the imaginary axis, the range of $s$ helps to take limits near zero and infinity. Notice that this is consistent with the spectral definition of the operator itself. This identity implies $(-\Delta_g)^sf(x)$ equals
\[\frac{1}{2\pi i}\int_{\gamma}z^{s}\int_Mf(y)\chi(x,y)u_0(x,y)F_0^z(d(x,y))d\textrm{vol}_g(y)dz\]
plus lower order terms.

Each parametrix summand can be expressed as
\[\frac{1}{2\pi i}\int_{\gamma}\int_Mz^{s+\frac{n}{4}-\frac{\nu+1}{2}}K_{n/2-\nu-1}(\sqrt{z}d_g(x,y))\chi u_{\nu}(x,y)f(y)d_g(x,y)^{-\frac{n}{2}+\nu+1}d\textrm{vol}_g(y)dz.\]
Making the change of variables $w=\sqrt{z}d_g(x,y)$ this double integral might be expressed as the product of a contour integral and a space integral. Notice that the contour integral is a constant not depending on the geometry, which coincides with the explicit one corresponding to the tori case. As a consequence the above double integral, up to a constant, is given by
\[\int_Mf(y)d(x,y)^{-n+2s+2\nu}\chi(x,y)u_{\nu}(x,y)d\textrm{vol}_g(y).\]
The other terms are less singular and can be handled similarly. Let us now turn to the error term which has the form
\[\frac{1}{2\pi i}\int_{\gamma}z^s(-\Delta_g-z)^{-1}R^zf(x)dz.\]
Subtitution of the expressions above leads to
\[\frac{1}{2\pi i}\int_{\gamma}z^s(-\Delta_g-z)^{-1}\int_M \chi(x,y)h_N(x,y)F^z_{N}(d_g(x,y))f(y)d\textrm{vol}_g(y)dz.\]
Let us denote by $H_0^s\subseteq H^s$ the subspace of functions orthogonal to constants. The claim follows combining the fact that $(-\Delta_g-z)^{-1}:H^s_0(M)\rightarrow H^s_0(M)$ with constant $O(1+(1+|z|)^{-1})$, Minkowski inequality, Sobolev's embedding $H^{n/2+\epsilon}\hookrightarrow L^{\infty}$, the recursive relations satisfied by the Bessel potentials and the uniform estimates for the modified Bessel function. For general $s\in (0,1)$ we apply  the above formula to $-\Delta_gf$, integrate by parts the principal part including the rest in the kernel $k_N$.

At this point a comment is in order: the Sobolev embedding theorem mentioned above can be proved independently of our next section. Indeed
\[f(x)=(-\Delta_g)^{-n/2-\epsilon}(-\Delta_g)^{n/2+\epsilon}f(x)=\sum_{\nu}\frac{a_{\nu}}{\lambda_{\nu}^{n/2+\epsilon}}Y_{\nu}(x)\]
where $(-\Delta_g)^{n/2+\epsilon}f(x)=\sum_{\nu}a_{\nu}Y_{\nu}(x)$ is the eingenfunction decomposition with $-\Delta_gY_k=\lambda_kY_k$. One may now apply Cauchy-Schwartz inequality, Weyl's law estimates and Plancherel to conclude $|f(x)|\leq C\|f\|_{H^{n/2+\epsilon}}$.

\remark[Manifolds with boundary] \normalfont In this case one uses Dirichlet or Neumann eigenfunctions and define acordingly its fractional operator. The parametrix still works due to its local character but the cut off $\chi$ should be taken more carefully. It would be enough to take $\chi(x,y)$ to be supported in a ball inside the manifold, but as a consequence the bounds of its derivatives which appear in the constant within the error term degenerate as we approach the boundary. We leave the details to the reader.

\remark[The Sobolev embedding theorem] \normalfont It is well known for integral number of derivatives and, therefore, we may restrict ourselves without loss of generality to fractions of the Laplace-Beltrami operator  $(-\Delta_g)^s$ where $s\in (0,1)$. The proof follows the usual lines (cf. Brascamp and Lieb treatise \cite{BL}). Indeed, for any $f$ with zero mean
\[\|f\|_{L^p(M)}=\sup_{\|g\|_q=1}\left|\int_Mf(y)g(y)dy\right|=\sup_{\|g\|_q=1}\left|\int \Lambda^sf\Lambda^{-s}g\right|\leq\|f\|_{H^s(M)}\|\Lambda^{-s}g\|_{2}\]
where $p^{-1}+q^{-1}=1$ are conjugates. But the last $L^2$-norm is bounded since it equals $\int g(-\Delta_g)^{-s}g$ which, due to our formula, is susceptible to an application of the Hardy-Littlewood-Sobolev inequality. The error term introduced is even nicer. To show it does not affect the validity of our statement  it is enough to interpolate between the $L^2\rightarrow L^2$ and the $L^1\rightarrow L^{\infty}$ bounds. The former has already been settled, let us show how to deal with the latter for which the following estimate holds
\[\begin{split}
\|(-\Delta_g-z)^{-1}Ef\|_{L^{\infty}(M)}&\leq C\|(-\Delta_g-z)^{-1}Ef\|_{H^{n/2+\epsilon}(M)}\\
&\leq  C\|Ef\|_{H^{n/2+\epsilon}(M)}\leq C\|f\|_{L^1(M)}
\end{split}\]
where we are denoting by $Ef$ the space integral in the error term, whose kernel is able to absorb derivatives without affecting integrability.

\section{\textbf{Acknowledgments}}

The authors were partially supported by ICMAT-Severo Ochoa project SEV-2011-0087 and the MTM2011-2281 project of the MCINN (Spain).


\begin{thebibliography}{10}
\bibitem{AOCM}
Alonso-Or\'an, D.; C\'ordoba, A.; Mart\'inez, A. D., {\em Continuity of weak solutions of the critical quasigeostrophic equations on the sphere}, unpublished.%\bibitem{Au}
\bibitem{AOCM2}
Alonso-Or\'an, D.; C\'ordoba, A.; Mart\'inez, A. D., {\em Global well--posedness of critical surface quasigeostrophic equation on the sphere}, unpublished.
\bibitem{Au}
Aubin, T., {\em Nonlinear Analysis on Manifolds. Monge-Amp\`ere equations}, Springer 1982.
\bibitem{Au2}
Aubin, T., {\em Espaces de Sobolev sur les vari\'et\'es riemannienes}, Bull. Sci. Math (2) 100 (1976), pp. 149-173.
%\bibitem{BC}
%Balodis, P.; C\'ordoba, A., {\em An inequality for Riesz transforms implying blow-up for some nonlinear and nonlocal transport equations}, Adv. of Math. Vol. 214 (10) (2007), pp. 1-39.
\bibitem{BL}
Loss, M.; Lieb, E. H., {\em Analysis, second edition}, Graduate Studies in Mathematics 14, 1997.
%\bibitem{CSil} Caffarelli, L. A.; Silvestre, L., {\em An extension problem related to the fractional Laplacian}, Communications in Partial Differential Equations, 32 (2007) Issue 8, 1245.
\bibitem{CS}
Caffarelli, L. A.; Sire, Y., {\em On some pointwise inequalities involving nonlocal operators}, Arxiv 1604.05665.
\bibitem{CV} 
Caffarelli, L. A.; Vasseur, A., {\em Drift diffusion equations with fractional diffusion and the quasi-geostrophic equation}, Ann. of Math. Vol. 171, No. 3 (2010), pp. 1903-1930.
%\bibitem{CV2}
%Caffarelli, L. A.; Vasseur, A., {\em The De Giorgi Method for Nonlocal Fluid Dynamics}.
\bibitem{Ch}
Chavel, I., {\em Eigenvalues in Riemannian geometry}, Academic Press, 1984.
%\bibitem{Ch2}
%Chavel, I, {\em Riemannian geometry: a modern introduction}, Cambridge University Press, 1993.
%\bibitem{CL}
%Coddington, E. A.; Levinson, N., {\em Theory of Ordinary and Differential Equations}, McGraw-Hill, 1955.
\bibitem{CI1}
Constantin, P.; Ignatova, M., {\em Critical SQG in bounded domains}, Annals of PDE 2:8 (2016), 42 pp.
\bibitem{CI2}
Constantin, P.; Ignatova, M., {\em Nonlinear lower bounds for the fractional Laplacian with Dirichlet boundary conditions and applications},  to appear in Int. Math. Res. Notices, 2016.
%\bibitem{CMT}
%Constantin, P.; Majda, A.; Tabak, E., {\em Formation of strong fronts in the 2D quasi-geostrophic thermal active scalar}, Nonlinearity 7 (1994), pp. 1495-1533.
\bibitem{CVi}
Constantin, P.; Vicol, V., {\em Nonlinear maximum principles for dissipative linear nonlocal operators and applications}, Geometric and Functional Analysis, 22 (2012), no. 5, 1289-1321. 
%\bibitem{Cdie}
%C\'ordoba, D., {\em Nonexistence of simple hyperbolic blow-up for the quasi-geostrophic equation}, Ann. of Math. 148 (1998), pp. 1135-1152.
\bibitem{CC} 
C\'ordoba, A.; C\'ordoba, D., {\em A Maximum Principle Applied to Quasi-Geostrophic Equations}, Commun. Math. Phys. 249 (2004), pp. 511-528.
\bibitem{CC2} 
C\'ordoba, A.; C\'ordoba, D., {\em A pointwise estimate for fractionary derivatives with applications to partial differential equations.} Proceedings of the National Academy of Sciences of the United States of America 100 (26) (2003), pp. 15316-15317.
%\bibitem{CCF}
%C\'ordoba, A.; C\'ordoba, D.; Fontelos, M. A., {\em Formation of singularities for a transport equation with nonlocal velocity}, Ann. of Math., 162 (2005), pp. 1377-1389.
%\bibitem{CCG} 
%C\'ordoba, A.; C\'ordoba, D.; Gancedo, F., {\em Interface evolution: the Hele-Shaw and Muskat problems}, Ann. of Math. Vol. 173 (2011), pp. 477-542.
\bibitem{CM}
C\'ordoba, A.; Mart\'inez, A. D., {\em A pointwise inequality for fractional laplacians}, Adv. of Math., Vol. 280 (2015), pp. 79-85.
%\bibitem{DG}
%De Giorgi, E., {\em Sulla differenziabilit\`a e l'analiticit\`a delle estremali degli integrali multipli regolari}, Mem. Accad. Sci. Torino Cl. Sci. Fis. Mat. Nat. (3) 3 (1957), pp. 25-43.
\bibitem{DKS}
Dos Santos Ferreira, D.; Kenig, C.; Salo, M., {\em On $L^p$ resolvent estimates for Laplace--Beltrami operators on compact manifolds}, Forum Math. 26, pp. 815-849 (2014).
%\bibitem{DZ}
%Duandikoetxea, J., {\em An\'alisis de Fourier}, Ediciones UAM, 1989.
%\bibitem{LE}
%Evans, L. C., {\em Partial differential equations}, Graduate Studies in Mathematics Vol. 19, 1998.
%\bibitem{Fo}
%Folland, G. B., {\em Introduction to Partial Differential Equations}, second edition. Princeton University Press, 1995.
\bibitem{Ho}
H\"ormander, L., {\em The analysis of linear partial differential operators (vol. III)}, Grundlehren der mathematischen Wissenschaften, Vol. 274, Springer, 1985.
%\bibitem{Ke}
%Kellogg, O. D., {\em Foundations of Potential Theory}, Dover, 1929.
\bibitem{KNV}
Kiselev, A.; Nazarov, F.; Volberg, A., {\em Global well-posedness for the critical 2D dissipative quasigeostrophic equation}, Invent. Math. 167 (2007), pp. 445-453.
%\bibitem{JL}
%Lions, J.L., {\em Quelque M\'ethodes de R\'esolutions des Probl\`emes aux Limites Non-Lin\'eaires}, Dunod, Paris, 1969.
\bibitem{LY}
Li, P; S-T. Yau, {\em On the parabolic kernel of the Schr\"odinger operator}, Acta Math. n.3-4, pp. 153-201 (1986).
\bibitem{MP}
Minakshisundaram, S.; Pleijel, A., {\em Some properties of the eigenfunctions of the Laplace-operator on Riemannian manifolds}, Canadian J. Math. 1 (1949), 242-256.
%\bibitem{N}
%Nash, J., {\em Continuity of Solutions of Parabolic and Elliptic Equations}, American Journal of Mathematics
%Vol. 80, No. 4 (1958), pp. 931-954.
%\bibitem{PW}
%Protter, M. H.; Weinberger, H. F., {\em Maximum principles in differential equations}, Pretice-Hall, 1967.
%\bibitem{R}
%Resnick, S. G. {\em Dynamical problems in non-linear advective partial differential equations}, PhD thesis (Chicago University, 1995).
%\bibitem{So}
%Sogge, Ch. D., {\em The Hangzhou lectures on Eigenfunctions of the Laplacian}, Princeton University Press, 2014.
%\bibitem{S0}
%Stein, E. M., {\em Harmonic Analysis: Real-Variable Methods, Orthogonality and Oscillatory Integrals}, Princeton University Press, 1993.
\bibitem{S1}
Stein, E. M., {\em Singular integrals and Differentiability Properties of Functions}, Princeton University Press, 1970.
%\bibitem{St}
%Strichartz, R. S., {\em Analysis of the Laplacian on the complete Riemannian Manifolds}, Journal of Functional Analysis 52 (1983), pp. 48-79.
%\bibitem{V}
%Vasseur, A. F., {\em The De Giorgi method for elliptic and parabolic equations and some applications}, to appear in Lectures on the Analysis of Nonlinear Partial Differential Equations Vol. 4. 
%\bibitem{W}
%Widman, K., {\em Inequalities for the Green function and boundary continuity of the gradient of solutions of elliptic differential equations}, Math. Scan. 21 (1967), pp. 17-37.
%\bibitem{SS} 
%Seeger, A.; Sogge, C. D., {\em Bounds for eigenfunctions of differential operators}, Indiana U. Math. J., Vol. 38 (1989), No. 38.
\bibitem{W} Watson, G. N., {\em A treatise on the theory of Bessel functions}, Cambridge University Press (1922).
\end{thebibliography}
\end{document}